\documentclass[11pt,amssym,twoside]{article}
\usepackage{amssymb}
\newtheorem{thm}{Theorem}[section]

\newtheorem{lem}[thm]{Lemma}

\newtheorem{prop}[thm]{Proposition}

\newtheorem{conj}[thm]{Conjecture}

\newcommand{\Cplx}{\mathbb C}

\begin{document}

\title{A geometric formulation of Siegel's diophantine theorem}%
\author{Arash Rastegar}%


\maketitle
\begin{abstract}
In this paper, we introduce an algebro-geometric formulation for Siegel's theorem 
using an improvement of Lang's version of Roth's theorem over finitely generated fields of 
characteristic zero [Lan]. In fact, we prove that, for an affine open curve $X$
in an irreducible smooth curve of genus $\geq 1$, any finitely generated subgroup of the additive
group of the affine ambient space intersects $X$ in only finitely many points.
This was proved only for finitely generated subgroups defined over a localization of the 
ring of integers of a number field by Mahler and others.
\end{abstract}
\section*{Introduction}
Siegel [Sie] has shown that an affine curve with coefficients in a 
number field and of genus $\geq 1$ has only a finite number of points whose coordinates 
are integers of that field. Mahler has conjectured that a similar statement holds 
for points having only a finite number of primes in their denominators, and proved this 
for curves of genus one over the rationals by his p-adic analogue of the Thue-Siegel theorem in [Mah1] and [Mah2]. 
Mahler's conjecture is fully proved today.

Mahler was a student in Frankfurt (1923-25) and Gottingen (1925-33) when he learned from 
Siegel about Thue's theorem and its improvements and generalisations.
Emmy Noether introduced him to the theory of p-adic numbers. He combined these two ideas in 1931 when he found an analogue of the 
Thue-Siegel theorem that involved both real and p-adic algebraic numbers. 
In 1955, Roth obtained his theorem on the rational approximations
of a real algebraic number. It was immediately clear to Mahler that his method should
also work for p-adic algebraic numbers. Some interesting work of this kind was
in fact carried out by D. Ridout, a student of Roth.

The method of Thue-Siegel-Roth has one fundamental disadvantage, that
of its non-effectiveness. The proof is entirely non-constructive, and by its
very nature does not lead to any upper bounds for possible solutions. Only
in some very special cases effective methods are known. Some are due
to Skolem and Gelfond.

In view of Roth's result, and the progress which had been made in the 
theory of abelian varieties (especially the Jacobian) since Siegel and Mahler's papers 
appeared, it seemed to Lang worth while to reconsider the question, which automatically carries with it a proof 
of Mahler's conjecture. The Jacobian is used in order to take a pull-back over the 
given curve of the standard covering given by $u\mapsto mu+a$ where m is a large integer, 
and $a$ is a suitable translation. 
Aside from Roth's theorem, Lang used only 
the classical properties of heights and the weak Mordell-Weil theorem. 
Mordell's conjecture which was proved now by Faltings, and then by Vojta and Bombieri,
claiming that a curve of genus $\geq 2$ has only a finite 
number of rational points would of course supersede the Siegel-Mahler theorem for such 
curves, but Lang conjectured that the latter holds in fact for abelian varieties: If $A$ is 
an abelian variety defined over a number field $K$, if $U$ is an open affine subset, and $R$ a 
subring of $K$ of finite type over $\mathbb{Z}$, then there is only a finite number of points of $U$ in $R$. 
The difficulty in trying to extend the proof to abelian varieties lies in the fact that there 
is a whole divisor at infinity, whereas for curves, there is only a finite number of points, 
which are all algebraic. This prediction by Lang, shows that he expected the Siegel's theorem be an 
algebro-geometric fact rather than an arithmetic proposition. 

Lang [Lan] even proved a version of Roth's theorem in case of finitely generated fields of characteristic zero,
but this is not strong enough to prove the geometric version of Siegel's theorem he had in mind.
Lang used the theory of heights which he and Neron had already generated for finitely generated fields of characteristic zero [Lan-Ner].

In this paper, we prove the stronger version of Roth's theorem that Lang needed
to improve Siegel's theorem. The idea is that, the maps of the form
$u\mapsto mu+a$ on an abelian variety are height increasing if $m$ is an integer $\geq 2$. 
This works even for height defined by Lang and Neron for finitely generated fields of characteristic zero.
We also use covering of  
a finitely generated group by images of finitely many maps of the above form.
As a reward, we get a geometric improvement of Siegel's theorem. Here's the statement of our generalization:

\begin{thm}
Let $X$ be an affine open subcurve of a connected smooth projectuve curve of genus $\geq 1$ defined
over $\mathbb{C}$  in the ambient affine space $\mathbb{A}^n(\mathbb{C})$ and let $F\subset \mathbb{A}^n(\mathbb{C})$
denote any finitely generated subgroup of $\mathbb{C}^n$. Then $X(K )\cap F$ is finite.
\end{thm}

This implies that Siegel's theorem is an algebro-geometric fact, not an arithmetic one.
Lang had the same geometric expectation when he formulated his conjecture that a curve of genus $\geq 2$ in its Jacobian should
intersect any finitely generated subgroup of the Jacobian in a finite set.
He even conjectured that divisible group of any finitely generated subgroup of the Jacobian intersects the embedded curve 
in finitely many points. Getting such a result is beyound our reach. Since we have only access to a result for finitely generated subgroups
defined over a finitely generated field. Therefore we state a conjecture following geometric philosophy of Lang as follows:

\begin{conj}
Let $X$ be an affine open subcurve of a connected smooth projectuve curve of genus $\geq 2$ defined
over $\mathbb{C}$  in the ambient affine space $\mathbb{A}^n(\mathbb{C})$ and let $F\subset \mathbb{A}^n(\mathbb{C})$
denote any finitely generated subgroup of $\mathbb{C}^n$ and $Div(F)$ denote the divisible subgroup of $\mathbb{C}^n$ associated to $F$. 
Then $X(K )\cap Div(F)$ is finite.
\end{conj}

In case $F$ is defined over a number field, this is proved by Mahler [Mah3] for any algebraic curve as above defined over $\mathbb{C}$.
Proving this conjecture would need a geometrization of Mahler's ideas over a finitely generated fields of characteristic zero.
\section{Diophantine approximation by subgroups of $\mathbb{C}^n$}

This section is devoted to proving theorems which were mentioned
in the introduction. The arguments are along the same lines as
analogous classical results.

Roth's theorem on Diophantine approximation of rational points on
projective line implies a version on projective varieties defined
over number-fields. 

\begin{thm} (Improvement of Roth's theorem on diopphantine approximation)
Fix a finitely generated field of characteristic zero $K$ 
and $\sigma :K\hookrightarrow \Cplx$ a
complex embedding. Let $A$ be an abelian variety
defined over $K$ and let $L$ be an very ample line-bundle on
$A$. Denote the arithmetic height function associated to the
line-bundle $L$ by $h_L$. Suppose $F\subset A(K)$ is a finitely generated subgroup.
Fix a Riemannian metric on $A_{\sigma}(\Cplx)$ and
let $d_{\sigma}$ denote the induced metric on
$A_{\sigma}(\Cplx)$. 
Then, for every $\delta>0$ and every choice
of an algebraic point $\alpha\in A(\bar {K})$ 
and all choices of a
constant $C$, there are only finitely many points
$\omega\in F$ approximating $\alpha$ such that 
$$
d_{\sigma}(\alpha ,\omega)\leq Ce^{-\delta h_L(\omega)}.
$$
\end{thm}

\begin{prop}
With assumptions of the above theorem, suppose for some $\delta_0>0$ 
we have that, for any choice of a constant $C$
and every choice
of an algebraic point $\alpha\in A(\bar {K})$
there are only finitely many points
$\omega\in F$ approximating $\alpha$ in the following manner
$$
d_{\sigma}(\alpha ,\omega)\leq Ce^{-\delta_0 h_L(\omega)}.
$$
Then, for every $\delta>0$ and every choice
of an algebraic point $\alpha\in A(\bar {K})$ and all choices of a
constant $C$, there are only finitely many points
$\omega\in F$ approximating $\alpha$ such that 
$$
d_{\sigma}(\alpha ,\omega)\leq Ce^{-\delta h_L(\omega)}.
$$
\end{prop}
\textbf{Proof (Proposition).} Note that, we have assumed that the above is true for some
$\delta_0>0$ 
without any assumption  
on $\alpha$.
Let $\delta'>0$ be infimum of such $\delta_0>0$.
The subset $F$ is disjoint union of the images of finitely many height-increasing
self-endomorphisms $\phi_i:A(K)\to A(K)$ 
defined over $K$ such that for
all $i$ we have
$$
h_L(\phi_i(f))=mh_L(f)+O(1)
$$
where $m>1$. Take $\phi_i:A(K)\to A(K)$ to be of the form $u \mapsto mu+a_i$
where $a_i$ are representatives of the finite group quotient $F/mF$.

Fix $\epsilon>0$ such that $\epsilon<\delta' <m\epsilon$. 
Suppose that $w_n$ is an infinite sequence of elements in
$F$ such that $\omega_n\to \alpha$ which satisfies the estimate
$$
d_{\sigma}(\alpha ,\omega_n)\leq Ce^{-\epsilon h_L(\omega_n)}.
$$
then infinitely many of them are images of elements of $F$ under
the same $\phi_i$. By going to a subsequence, one can find a
sequence $\omega'_n$ in $F$ and an algebraic point $\alpha'$ in
$A(\bar {K})$ such that $\omega'_n \to \alpha'$ and for a fixed
$\phi_i$ we have $\phi_i(\alpha')=\alpha$ and
$\phi_i(\omega'_n)=\omega_n$ for all $n$. Then
$$
d_{\sigma}(\alpha ,\omega_n)\leq Ce^{-\epsilon h_L(\omega_n)}\leq
C'e^{-\epsilon m_i h_L(\omega'_n)}
$$
for an appropriate constant $C'$. On the other hand,
$$
d_{\sigma}(\alpha' ,\omega'_n)\leq C''d_{\sigma}(\alpha ,\omega_n)
$$
holds for an appropriate constant $C''$ and large $n$ by
injectivity of $d\phi_i^{-1}$ on the tangent space of $\alpha$.
This contradicts our assumption on $\delta'$, because $\delta' <m_i\epsilon$.
$\square$ 
\\
\textbf{Proof (Theorem).}  If we assume that points of $F$ and covering maps are defined over some
number-field, Roth's theorem implies that the assumption of theorem is true for any 
$\delta_0>2$. The same is true for finitely generated fields of characteristic zero
by a result of Lang [Lan] generalizing Roth's theorem.$\square$

\section{Geometric formulation of Siegel's theorem}

Let us state a more precise version of our version of Siegel's theorem.

\begin{thm}(Geometric version of Siegel's theorem on integral points) Fix a
finitely generated field of characteristic zero $K$. 
Let $X$ be an affine open subcurve of a connected smooth projectuve curve of genus $\geq 1$ defined
over $K$  in the ambient affine space $\mathbb{A}^n(K)$ 
and let $F\subset \mathbb{A}^n(K)$
denote any finitely generated subgroup of $K^n$. 
Then $X(K )\cap F$ is finite.
\end{thm}

We borrow a classical lemma [Ser] whose proof goes very similar to that reference.

\begin{lem}
Let $K$ be a finitely generated field of characteristic zero. Let $X$ be a curve defined over $K$.
Assume genus of $X$ is $\geq 1$.
If $P_n$ is a sequence of distince
points in $X(K)$, which means that their heights tends to infinity and if
$\phi$ defined over $K$ is a non-constant rational function on $X$. From some point
on, no $P_n$ is pole of $\phi$. 
Then for $z_n=\phi(P_n)$ which a point of the
projective space defined over $K$ we have
$$
\lim_{n\to \infty} {{log|z_n|_{v}}\over {log H(z_n)}}=0
$$
\end{lem}
\textbf {Proof.}
Assume this is false. By taking a subsequence and replacing $\phi$ by
$1/\phi$, we may suppose that
$$
{{log|z_n|_{v}} \over {log H(z_n)}}\rightarrow \lambda
$$
where $-\infty <\lambda< 0$.
In particular, $z_n \rightarrow 0$ in $K_v$ and by taking a subsequence, we may assume
that $P_n$ converges to a zero $P_0$ of $\phi$. As we are on a curve, $P_0$ is an algebraic
point of X.
Between $H(P_n )$, the height corresponding to a morphism $X \rightarrow \mathbb P_N$,
and $H(z_n )$, corresponding to a morphism $X \rightarrow \mathbb P_1$, we have an inequality
$$
H(z_n )\ll H(P_n )^{l}
$$
for some positive $l$.
On the other hand, if $e$ is the multiplicity of $P_0$ as a zero of $\phi$, we have
$|z_n|_{v}\approx d_{v}(P_n, P_0)^{e}$.
Therefore, there is $c > 0$ such that for sufficiently large $n$,
$$
d_{v}(P_n, P_0)\leq 1/H(P_n )^c
$$
which contradicts the approximation theorem.$\square$
\\
\textbf {Proof (Theorem).} Let $\sigma :K\hookrightarrow \Cplx$
denote a complex embedding of $K$. Fix a Riemannian metric on
$Jac(X)_{\sigma}(\Cplx)$ and let $d_{\sigma}$ denote the induced
metric on $Jac(X)_{\sigma}(\Cplx)$. Now, embed $X$ in $Jac(X)$. Then, by our version of Roth's
theorem, for every $\delta>0$ and every choice of an algebraic
point $\alpha\in X(\bar {K})$ and all choices of a constant $C$, there
are only finitely many points $\omega\in F\cap X(K)$ approximating
$\alpha$ such that
$$
d_{\sigma}(\alpha ,\omega)\leq H_L(\omega)^{-\delta }.
$$
where $log(H_L)=h_L$. In case $K$ is trancendental, we have to pick a model for $K$ over algebraic closure of 
$\mathbb{Q}$ in $K$ following Lang [Lan].

Now if $P_n$ is a sequence of distince
points in $X(K)\cap F$, their heights tends to infinity and if
$\phi$ is a non-constant rational function on $X$ from some point
on no $P_n$ is pole of $\phi$.
Then by above lemma
$$
\lim_{n\to \infty} {{log|z_n|_{\sigma}}\over {log H(z_n)}}=0
$$
On the other hand, one defines height of points defined over $K$ by
$$
H(z)=\prod_{v\in M_K} sup(1,|z|_v),
$$
where $|.|_v$ are normalized according to a product formula.
Since covering maps of $F$ are height expanding, we know 
that $F$ is forward orbit of finitely many points. So for a
finite set of places $S$ we have
$$
H(z)=\prod_{v\in S} sup(1,|z|_v),
$$
and therefore
$$
log H(z)=\sum_{v\in S} log(sup(1,|z|_v)).
$$
Then, we have 
$$
1=\sum_{v\in S} sup(0,{{log|z_n|_{\sigma}}\over {log H(z_n)}})\leq \sum_{v\in S}  {{log|z_n|_{\sigma}}\over {log H(z_n)}}
$$
which could not be true, because the above limit is zero. This
implies the finiteness result we are seeking for. $\square$

\subsection*{acknowledgements}
I have benefited from conversations with M. Hadian, A. Rajaei, P. Sarnak, N. Talebizadeh, 
for which I am thankful. Peter Sarnak particularly gave crucial comments which led 
to the final version of the paper. 
I would also like to thank Sharif University
of Technology for finantial support and Princeton University for warm hospitality.

Sharif University of Technology, e-mail: rastegar@sharif.edu
\\Princeton University, e-mail:
rastegar@princeton.edu

\end{document}